\documentclass[a4 paper,10pt,twoside]{article}
\pdfoutput=1
\usepackage[centertags]{amsmath}
\usepackage{amssymb}
\usepackage{amsfonts}
\usepackage{amsthm}
\usepackage{newlfont}
\usepackage{amscd}
\usepackage{tikz}
\usepackage[top=1.7in, bottom=1.7in, left=1.5in, right=1.5in]{geometry}
\title{\large On the reciprocal sum of Jacobsthal numbers}
\author{}
\date{}
\theoremstyle{plain}
\newtheorem{theorem}{Theorem}[section]

\newtheorem{lemma}[theorem]{Lemma}
\newtheorem{corollary}[theorem]{Corollary}

\theoremstyle{example}

\numberwithin{equation}{section}
\theoremstyle{definition}

\theoremstyle{remark}

\begin {document}
\maketitle
\begin{center}
\textbf{ Ahmed Gaber \vspace{0.15 cm}\\ Faculty of Science, Ain Shams
University, Egypt\\
a.gaber@sci.asu.edu.eg}
\end{center}
\vspace{0.1 cm}
\begin{abstract}
In this paper, we  study the reciprocal sums of
the Jacobsthal numbers. We establish many results on the infinite sum and alternating infinite sum
of the reciprocals of  Jacobsthal numbers  and square Jacobsthal numbers.
\end{abstract}
\begin{flushleft}
\end{flushleft}
\begin{flushleft}

\end{flushleft}
\section{Introduction}
 The Jacobsthal polynomials $J_{n}(x)$ were first studied by E. E. Jacobsthal around 1919.  Jacobsthal polynomials $J_{n}(x)$ can be defined by the  recurrence $J_{0}(x) = 0, J_{1}(x) = 1$ and $J_{n}(x) = J_{n-1}(x) + xJ_{n-2}(x)$, for $n \geq 2$. Clearly, $J_{n}(1) = F_{n}$, the sequence of Fibonacci numbers. The sequence $ J_{n}= J_{n}(2) =J_{n+2}=J_{n+1}+2J_{n}$ is the sequence of Jacobsthal numbers. The first few terms are $ \{0, 1,1, 3, 5, 11, 21, 43, 85, ...\}$. The first apperance of this sequence was in [4]. It was Horadam who first considered such a sequence in detail in his seminal paper [5] .  Motivated by Horadam’s work, a lot of research has been conducted. For some recent works, see [1-3,8], for example. Recently, there has been an increasing interest in studying the reciprocal sums
of the Fibonacci numbers ($J_{n}(1)$). For instance, see [9-11].\\

In this article we confine ourselves the infinite sum and alternating infinite sum
of the reciprocals of  Jacobsthal numbers  and square Jacobsthal numbers. We obtain many intersting results in this direction.

We first state several well known results on Jacobsthal numbers, which will be used
throughout the article. The detailed expositions can be found  in  [6,7].
\begin{lemma}
 For $n\geq 1$, we have
 \begin{equation}\label{*}
 J_{n}+J_{n+1}=2^{n}
 \end{equation}
\end{lemma}
\begin{lemma}
 For any positive integer $n\geq 1$, we have
 \begin{equation}\label{*}
 J_{n}< 2^{n}
 \end{equation}
 For any positive integer $n\geq 2$, we have
 \begin{equation}\label{*}
 J_{n}< 2^{n-1}
 \end{equation}
  For any positive integer $n\geq 3$, we have
 \begin{equation}\label{*}
 2^{n-2}<J_{n}< 2^{n-1}
 \end{equation}
\end{lemma}
\begin{lemma}
 For $n\geq 1, k\geq1, n\geq k$, the Jacobsthal numbers has the Cassini-like formula
 \begin{equation}\label{*}
 J_{n+k}J_{n-k}-J_{n}^{2}= (-1)^{n-k+1}2^{n-k}J_{k}^{2}
 \end{equation}
\end{lemma}
\begin{lemma}
 For $n\geq 1$, we have
 \begin{equation}\label{*}
 J_{n+1}^{2}-J_{n}^{2} =2^{n+1}J_{n-1}
 \end{equation}
\end{lemma}
\begin{lemma}
 For $n\geq 1$, we have
 \begin{equation}\label{*}
J_{n+1}^{2}+2J_{n}^{2}= J_{2n+1}
 \end{equation}
\end{lemma}

\section{Reciprocal sums of the Jacobsthal numbers }
 In this section,

\begin{theorem}
    Let $n\geq 2$. Then \\
   $$ J_{n-2} < (\displaystyle \sum_{k=n}^{\infty} \frac{1}{J_{k}})^{-1}< 4(J_{n-2}+1).$$
\end{theorem}

\begin{proof}
For $k\geq 1$, we have
\begin{eqnarray*}
\frac{1}{J_{n}}- \frac{2}{J_{n+2}}-\frac{1}{J_{n+3}}&=&\frac{J_{n+2}-2J_{n}}{J_{n}J_{n+2}}-\frac{1}{J_{n+3}}\\
&=& \frac{J_{n+1}}{J_{n}J_{n+2}}-\frac{1}{J_{n+3}}\\
&=&\frac{J_{n+1}J_{n+3}- J_{n}J_{n+2}}{J_{n}J_{n+2}J_{n+3}}
\end{eqnarray*}
Clearly, $J_{n+1}J_{n+3}- J_{n}J_{n+2}> 0$. Therefore,\\
$$\frac{1}{J_{n}}> \frac{1}{J_{n+2}}+ \frac{1}{J_{n+2}} +\frac{1}{J_{n+3}}.$$
Repeating this inequality for $n> 2$, we obtain
\begin{eqnarray*}
\frac{1}{J_{n-2}} &>& \frac{1}{J_{n}}+ \frac{1}{J_{n}} +\frac{1}{J_{n+1}}\\
&>& \frac{1}{J_{n}}+ \frac{1}{J_{n+1}} +(\frac{1}{J_{n+2}}+ \frac{1}{J_{n+2}}+\frac{1}{J_{n+3}})\\
&>& \frac{1}{J_{n}}+ \frac{1}{J_{n+1}} + \frac{1}{J_{n+2}}+\frac{1}{J_{n+3}}+(\frac{1}{J_{n+4}}+ \frac{1}{J_{n+4}}+\frac{1}{J_{n+5}})\\
&>& ...\\
&>& \displaystyle \sum_{k=n}^{\infty} \frac{1}{J_{k}}.
\end{eqnarray*}
Therefore, for $n > 2$, we have
\begin{equation}\label{*}
 \displaystyle \sum_{k=n}^{\infty} \frac{1}{J_{k}}< \frac{1}{J_{n-2}}.
\end{equation}
Now assume that $k\geq m\geq 1$. Then
 \begin{center}
 $ \frac{J_{k-m}+1}{J_{k}}\geq \frac{2^{k-m-2}+1}{2^{k-1}}> 2^{-m-1}.$
 \end{center}
Let $m=k-n+2$. Then
 \begin{equation*}
  \displaystyle \sum_{k=n}^{\infty} \frac{J_{n-2}+1}{J_{k}} > \displaystyle \sum_{k=n}^{\infty}2^{n-k-3}= \displaystyle \sum_{t=2}^{\infty}2^{-t-1}= \frac{1}{4}.
\end{equation*}
Therefore,
\begin{equation}\label{*}
  \displaystyle \sum_{k=n}^{\infty} \frac{1}{J_{k}} >\frac{1}{4(J_{n-2}+1)}.
\end{equation}
Combining (2.1) and (2.2), we get
\begin{equation*}
\frac{1}{4(J_{n-2}+1)} < \displaystyle \sum_{k=n}^{\infty} \frac{1}{J_{k}}< \frac{1}{J_{n-2}}.
\end{equation*}
Hence,
\begin{equation*}
J_{n-2} < (\displaystyle \sum_{k=n}^{\infty} \frac{1}{J_{k}})^{-1}< 4(J_{n-2}+1).
\end{equation*}
This completes the proof.
\end{proof}
 \begin{theorem}
 Let $n \in \mathbf{N}$. If $n$ is odd, then\\
  $$\lfloor (\displaystyle \sum_{k=n}^{\infty} \frac{1}{J_{k}^{2}})^{-1}\rfloor \leq  J_{n-1}J_{n}.$$
   \end{theorem}
\begin{proof}
 If $n=1$, then $\displaystyle \sum_{k=1}^{\infty} \frac{1}{J_{k}^{2}}> \frac{1}{J_{1}^{2}}=1$. This implies that $0< (\displaystyle \sum_{k=n}^{\infty} \frac{1}{J_{k}^{2}})^{-1}< 1$. Therefore, $\lfloor (\displaystyle \sum_{k=1}^{\infty} \frac{1}{J_{k}^{2}})^{-1}\rfloor = 0= J_{0}J_{1}.$ Let $n\geq 3$. Then
 \begin{eqnarray*}
\frac{1}{J_{n-1}J_{n}}- \frac{1}{J_{n}^{2}}-\frac{2}{J_{n+1}^{2}}-\frac{4}{J_{n+1}J_{n+2}} &=& \frac{J_{n}-J_{n-1}}{J_{n-1}J_{n}^{2}}- \frac{2[J_{n+2}+2J_{n+1}]}{J_{n+1}^{2}J_{n+2}} \\
&=& \frac{2J_{n-2}}{J_{n-1}J_{n}^{2}}- \frac{2J_{n+3}}{J_{n+1}^{2}J_{n+2}} \\
&=& \frac{2J_{n-2}J_{n+1}^{2}J_{n+2}-2J_{n-1}J_{n}^{2}J_{n+3}}{J_{n-1}J_{n}^{2}J_{n+1}^{2}J_{n+2}} \\
&=& \frac{2J_{n+1}^{2}[J_{n}^{2}+(-1)^{n-1}2^{n-2}J_{2}^{2}]-2J_{n}^{2}[J_{n+1}^{2}+(-1)^{n}2^{n-1}J_{2}^{2}]}{J_{n-1}J_{n}^{2}J_{n+1}^{2}J_{n+2}}\\
&=& \frac{(-1)^{n-1}2^{n-1}J_{n+1}^{2}+(-1)^{n+1}2^{n}J_{n}^{2}}{J_{n-1}J_{n}^{2}J_{n+1}^{2}J_{n+2}}\\
&=& \frac{(-1)^{n-1}2^{n-1}(J_{n+1}^{2}+2J_{n}^{2})}{J_{n-1}J_{n}^{2}J_{n+1}^{2}J_{n+2}}\\
&=& \frac{(-1)^{n-1}2^{n-1}J_{2n+1}}{J_{n-1}J_{n}^{2}J_{n+1}^{2}J_{n+2}}\\
\end{eqnarray*}
Since $n$ is odd, then
$$ \frac{1}{J_{n-1}J_{n}} > \frac{1}{J_{n}^{2}}+\frac{2}{J_{n+1}^{2}}+\frac{4}{J_{n+1}J_{n+2}}.$$
Repeating this for $n\geq 3$, we obtain
\begin{eqnarray*}
\frac{1}{J_{n-1}J_{n}} &>& \frac{1}{J_{n}^{2}}+\frac{2}{J_{n+1}^{2}}+\frac{4}{J_{n+1}J_{n+2}}\\
&>& \frac{1}{J_{n}^{2}}+\frac{2}{J_{n+1}^{2}}+ 4(\frac{1}{J_{n+2}^{2}}+\frac{2}{J_{n+3}^{2}}+\frac{4}{J_{n+3}J_{n+4}})\\
&>& \frac{1}{J_{n}^{2}}+\frac{2}{J_{n+1}^{2}}+ \frac{4}{J_{n+2}^{2}}+\frac{8}{J_{n+3}^{2}}+16 (\frac{1}{J_{n+4}^{2}}+\frac{2}{J_{n+5}^{2}}+\frac{4}{J_{n+5}J_{n+6}})\\
&>& ...\\
&>& \displaystyle \sum_{k=0}^{\infty} \frac{2^{k}}{J_{n+k}^{2}}\\
&>& \displaystyle \sum_{k=0}^{\infty} \frac{1}{J_{n+k}^{2}}\\
&=& \displaystyle \sum_{k=n}^{\infty} \frac{1}{J_{k}^{2}}
\end{eqnarray*}
Therefore,
\begin{equation*}
 \displaystyle \sum_{k=n}^{\infty} \frac{1}{J_{k}^{2}} <   \frac{1}{J_{n-1}J_{n}}
\end{equation*}
Hence, $\lfloor (\displaystyle \sum_{k=n}^{\infty} \frac{1}{J_{k}^{2}})^{-1}\rfloor \leq  J_{n-1}J_{n}.$
\end{proof}
\section{Alternating reciprocal sums of the Jacobsthal numbers }
\begin{theorem}
 Let $n \in \mathbf{N}$. If $n$ is even, then\\
  $$\lfloor (\displaystyle \sum_{k=n}^{\infty} \frac{(-1)^{k}}{J_{k}^{2}})^{-1}\rfloor =  2^{n-1}-1.$$
   \end{theorem}
\begin{proof}
   For $n \geq 1$, we have
   \begin{eqnarray*}
\frac{(-1)^{n}}{J_{n-1}+J_{n}-(-1)^{n}}-\frac{(-1)^{n}}{J_{n}} - \frac{(-1)^{n+1}}{J_{n}+J_{n+1}-(-1)^{n+1}} &=& \frac{(-1)^{n+1} J_{n-1}+1}{J_{n}(J_{n-1}+J_{n}-(-1)^{n})}+ \frac{(-1)^{n}}{J_{n}+J_{n+1}-(-1)^{n+1}} \\
&=& \frac{M}{J_{n}(J_{n-1}+J_{n}-(-1)^{n})(J_{n}+J_{n+1}-(-1)^{n+1})}
\end{eqnarray*}
Where $M= (-1)^{n+1}J_{n-1}J_{n+1}+J_{n+1}-J_{n-1}+(-1)^{n}J_{n}^{2}+(-1)^{n}$. We simplify $M$ as follows
\begin{eqnarray*}
M &=&(-1)^{n+1}J_{n-1}J_{n+1}+J_{n+1}-J_{n-1}+(-1)^{n}J_{n}^{2}+(-1)^{n} \\
&=& (-1)^{n+1}(J_{n}^{2}+(-1)^{n}2^{n-1})+2^{n-1}+(-1)^{n}J_{n}^{2}+(-1)^{n} \\
&=& (-1)^{n}
\end{eqnarray*}
Since $n$ is even, then
\begin{eqnarray*}
\frac{1}{(-1)^{n}(J_{n-1}+J_{n})-1} &> & \frac{(-1)^{n}}{J_{n}} + \frac{1}{(-1)^{n+1}(J_{n}+J_{n+1})-1}
\end{eqnarray*}
A repetition of this shows that,
\begin{equation}\label{*}
  \displaystyle \sum_{k=n}^{\infty} \frac{(-1)^{k}}{J_{k}} < \frac{1}{(-1)^{n}(J_{n-1}+J_{n})-1}
\end{equation}
On the other hand, for $n\geq 1$, we find
  \begin{eqnarray*}
\frac{(-1)^{n}}{J_{n}-(-1)^{n}}-\frac{(-1)^{n}}{J_{n-1}+J_{n}} + \frac{(-1)^{n+1}}{J_{n}+J_{n+1} } &=& (-1)^{n}[\frac{ J_{n-1}}{J_{n}(J_{n-1}+J_{n} )}- \frac{1}{J_{n}+J_{n+1}}] \\
&=& (-1)^{n}[\frac{ J_{n-1}J_{n+1}-J_{n}^{2}}{J_{n}(J_{n-1}+J_{n} )(J_{n}+J_{n+1})}] \\
&=& \frac{ 2^{n-1}}{J_{n}(J_{n-1}+J_{n} )(J_{n}+J_{n+1})} \\
&>& 0
\end{eqnarray*}
Therefore,
\begin{equation}\label{*}
  \displaystyle \sum_{k=n}^{\infty} \frac{(-1)^{k}}{J_{k}} > \frac{1}{(-1)^{n}(J_{n-1}+J_{n})}.
\end{equation}
Combining (3.1) and (3.2), we obtain ( for $n$ even)
\begin{equation*}
J_{n-1}+J_{n}-1 < (\displaystyle \sum_{k=n}^{\infty} \frac{(-1)^{k}}{J_{k}})^{-1} < J_{n-1}+J_{n}.
\end{equation*}
Hence,
\begin{equation*}
  \lfloor(\displaystyle \sum_{k=n}^{\infty} \frac{(-1)^{k}}{J_{k}})^{-1}\rfloor = J_{n-1}+J_{n}-1=2^{n-1}-1.
\end{equation*}
\end{proof}
\begin{corollary}
  If $n$ is odd, then
  $$\lfloor(\displaystyle \sum_{k=n}^{\infty} \frac{(-1)^{k}}{J_{k}})^{-1}\rfloor \leq -(2^{n-1}+1).$$
\end{corollary}
\begin{proof}
  If $n$ is odd, then $M=-1$. We conclude that
  $$\frac{1}{(-1)^{n}(J_{n-1}+J_{n})-1} < \frac{(-1)^{n}}{J_{n}} + \frac{1}{(-1)^{n+1}(J_{n}+J_{n+1})-1}$$
Applying this repeatedly and invoking the relation $J_{n-1}+J_{n}=2^{n-1}$ give the desired result.
\end{proof}
\begin{theorem}
Let $n$ be a positive integer. Then
$\lceil (\displaystyle \sum_{k=n}^{\infty} \frac{(-1)^{k}}{J_{k}^{2}})^{-1} \rceil \leq J_{n-1}^{2}+J_{n}^{2}-1.$
\end{theorem}
\begin{proof}
\begin{eqnarray*}
\frac{(-1)^{n}}{J_{n-1}^{2}+J_{n}^{2}-(-1)^{n}}-\frac{(-1)^{n}}{J_{n}^{2}} - \frac{(-1)^{n+1}}{J_{n}^{2}+J_{n+1}^{2}-(-1)^{n+1}} &=& \frac{(-1)^{n+1} J_{n-1}^{2}+1}{J_{n}^{2}(J_{n-1}^{2}+J_{n}^{2}-(-1)^{n})}+ \frac{(-1)^{n}}{J_{n}^{2}+J_{n+1}^{2}+(-1)^{n}} \\
&=& \frac{N}{J_{n}^{2}(J_{n-1}^{2}+J_{n}^{2}-(-1)^{n})(J_{n}^{2}+J_{n+1}^{2}+(-1)^{n})}
\end{eqnarray*}
Where $N= (-1)^{n+1}J_{n-1}^{2}J_{n+1}^{2}+J_{n+1}^{2}-J_{n-1}^{2}+(-1)^{n}J_{n}^{4}+(-1)^{n}.$ Now, we bound $N$.

 Using $J_{n-1}J_{n+1}= J_{n}^{2}+(-1)^{n}2^{n-1}$, we immediately get
 \begin{equation*}
 N= 2^{n+1}J_{n-1}+J_{n}^{2}+(-1)^{n+1}2^{2n-2}-J_{n-1}^{2}-2^{n}J_{n}^{2}+(-1)^{n}.
 \end{equation*}
 Elementary manipulations of the inequality (1.4) entail that
 \begin{eqnarray*}
N &<& 2^{2n-1}+ 2^{2n-2}+2^{2n-2}-2^{2n-6}-2^{3n-4}+1 \\
&=& \frac{63}{64} 2^{2n}-2^{3n-4}+1\\
&<& 2^{2n}-2^{3n-4}+1.
\end{eqnarray*}
We note that $2n < 3n-4$ for $n\geq 5$. Hence, $N< 0$ for $n\geq 5$. In other words, we have (for $n\geq 5$)
 \begin{equation*}
 \frac{1}{(-1)^{n}(J_{n-1}^{2}+J_{n}^{2})-1}<\frac{(-1)^{n}}{J_{n}^{2}} + \frac{1}{(-1)^{n+1}(J_{n}^{2}+J_{n+1}^{2})-1}.
 \end{equation*}
 using this equation, we can prove that
   \begin{equation*}
 \displaystyle \sum_{k=n}^{\infty} \frac{(-1)^{k}}{J_{k}^{2}}> \frac{1}{(-1)^{n}(J_{n-1}^{2}+J_{n}^{2}-1)}
  \end{equation*}
  For $n$ even, we get  $(\displaystyle \sum_{k=n}^{\infty} \frac{(-1)^{k}}{J_{k}^{2}})^{-1} < J_{n-1}^{2}+J_{n}^{2}-1$. The result follows.
\end{proof}

\end{document}